\documentclass{amsart}

\newtheorem{theorem}{Theorem}[section]
\newtheorem{proposition}[theorem]{Proposition}
\newtheorem{corollary}[theorem]{Corollary}
\newtheorem{definition}[theorem]{Definition}

\begin{document}
\title[Spectrum of complex banded matrices]{On the discrete spectrum
of complex banded matrices }

\author{Leonid Golinskii}
\address{Institute for Low Temperature Physics and Engineering\\47\\ Lenin
ave.\\ Kharkiv \\ Ukraine} \email{golinskii@ilt.kharkov.ua}

\author{Mikhail Kudryavtsev}
\address{Institute for Low Temperature Physics and Engineering\\47\\ Lenin
ave.\\ Kharkov \\ Ukraine} \email{kudryavtsev@ilt.kharkov.ua;
kudryavstev@onet.com.ua}

\thanks{The work of the first author was supported in part by
INTAS Research Network NeCCA 03-51-6637. }

\begin{abstract}
The discrete spectrum of complex banded matrices that are compact
perturbations of the standard banded matrix of order $p$ is under
consideration. The rate of stabilization for the matrix entries
sharp in the sense of order which provides finiteness of the
discrete spectrum is found. The $p$-banded matrix with the
discrete spectrum having exactly $p$ limit points on the interval
$(-2,2)$ is constructed. The results are applied to the study of
the discrete spectrum of asymptotically periodic Jacobi matrices.

\end{abstract}

\maketitle

 \section{Introduction.}\setcounter{section}{1}
\setcounter{equation}{0}

$ \ \  \ $     In the recent papers \cite{EG1,EG2} I. Egorova and
L. Golinskii studied the discrete spectrum of complex Jacobi
matrices such that the operators in $\ell^2( \mathbb{N})$,
$\mathbb{N}:= \{1,2,...\}$ generated by these matrices are compact
perturbations of the discrete laplacian. In turn, these papers are
the discrete version of the known Pavlov's theorems (\cite{P1,P2})
for the differential operators of second order on the semiaxis.
The sufficient conditions for the spectrum to be finite and empty,
the domains containing the discrete spectrum and the conditions
for the limit sets of the discrete spectrum were found. The goal
of this work is to extend these results to the case of operators,
generated by banded matrices.

Let us remind that an infinite matrix
$D=\|d_{ij}\|_{i,j=1}^\infty$ is called the {\it banded matrix of
order $p$ } or just $p$-banded if
\begin{equation}
\label{0.1} d_{ij}=0, \quad |i-j|>p, \qquad d_{ij}\ne 0, \quad
|i-j|=p, \qquad d_{ij} \in \mathbb{C}.
\end{equation}
According to this definition, the Jacobi matrices are banded
matrices of order $p=1$. Throughout the whole paper we assume that
\begin{equation}\label{0.2}
\lim_{i\to \infty}d_{i,i\pm p}=1, \qquad
\lim_{i\to\infty}d_{i,i\pm r}=0,\quad |r|<p,
\end{equation}
and so the operators in $\ell^2=\ell^2(\mathbb{N})$ generated by
matrices (\ref{0.1})--(\ref{0.2}) are compact perturbation of the
standard banded operator
\begin{equation}\label{0.11}
D_0:\qquad d_{i,i\pm p}=1, \qquad d_{ij}=0, \quad |i-j| \ne p,
\qquad D_0=S^p+(S^*)^p,
\end{equation}
where $S$ is the one-sided shift operator in $\ell^2$. It is well
known that the spectrum $\sigma(D_0)$ of $D_0$ is the closed
interval $[-2,2]$. According to the Weyl theorem (see, e.g.,
\cite{GK}) the spectrum of the perturbed operator $\sigma(D) =
[-2,2]\bigcup \sigma_{\rm d}(D)$, where the discrete spectrum
$\sigma_{\rm d}(D)$ is at most denumerable set of points of the
complex plane, which are eigenvalues of finite algebraic
multiplicity. All its accumulation points belong to the interval
$[-2,2]$. Let us denote by $E_D$ the limit set for the set
$\sigma_{\rm d}(D)$. So, $E_D=\emptyset$ means that the discrete
spectrum is finite.

Remind that the {\it convergence exponent} or {\it
Taylor-Besicovitch index} of a closed point set $F\subset [-2,2]$
is the value
$$\tau(F):=\inf\{\varepsilon >0:\ \sum_{j=1}^\infty
|l_j|^\varepsilon <\infty\},$$ where $\{l_j\}$ are the adjacent
intervals of $F$.

\begin{definition}\label{def1}
 We say that the matrix $D$ $(\ref{0.1})$ belongs to the class $
\mathcal{P}_p(\beta)$, $0<\beta<1$, if
\begin{equation}
\label{0.3} q_n:=\left|d_{n,n-p}-1\right|+\sum_{r={-p+1}}^{p-1}
\left| d_{n,n+r} \right| + \left|d_{n,n+p}- 1\right|\leq
C_1\,\exp(-C_2\, n^\beta)\,,
\end{equation}
$n\in\mathbb{N}$, with the constants $C_1,\, C_2\,>0$, depending
on $D$.
\end{definition}

The main result of the present paper is the following
\begin{theorem}
Let $D\in\mathcal{P}_p(\beta)$ where $0<\beta<\frac{1}{2}$. Then $E_D$
is a closed point set of the Lebesgue measure zero and its
convergence exponent satisfies
\begin{equation} \label{1.6} \dim E_D  \leq
\tau(E_D)\leq\frac{1-2\beta}{1-\beta}\,,
\end{equation}
where $\dim E_D$ is the Hausdorff dimension of $E_D$. Moreover, if
$D\in \mathcal{P}(\frac{1}{2})$ then $E_D=\emptyset$, i.e., the
discrete spectrum is finite.
\end{theorem}

It turns out that the exponent $1/2$ in theorem 1 is sharp in the
following sense.
\begin{theorem}
\label{theor2} For arbitrary $\varepsilon > 0$ and arbitrary
points $\nu_1, \nu_2, \ldots, \nu_p \in (-2,2)$ there exists an
operator $D\in \mathcal{P}_p(\frac{1}{2} - \varepsilon)$ such that
its discrete spectrum $\sigma_{\rm d}(J)$ is infinite and,
moreover,
$$E_D=\{\nu_1,\nu_2,\ldots\nu_p\}.$$
\end{theorem}

The theorems 1 and 2 are proved in Section 4, where the domains
which contain $\sigma_{\rm d}(D)$ are also found (under the weaker
assumptions than (\ref{0.3})). In Section 2 the connection is
established between the discrete spectrum and zeros of the
determinant constructed of $p$ linearly independent solutions of
the linear difference equations for the eigenvector. In Section 3
the properties of the matrix Jost solutions are studied. Finally,
in the last Sections 5 and 6 the main results are applied to study
the spectrum of doubly-infinite complex banded matrices and the
spectrum of doubly-infinite asymptotically $p$-periodic complex
Jacobi matrix.

\setcounter{section}{1} \setcounter{equation}{0}
\section {The determinants of independent solutions and the eigenvalues}

We start out with the equation
\begin{equation} \label{1.1} D \vec{y} = \lambda \vec{y}
\end{equation}
for generalized eigenvectors $\vec{y} = \{y_n\}_{n\ge 1}$ in the
coordinate form:
\begin{equation} \label{1.2}
\left\{\begin{array}{l}
d_{11}\,y_1 + d_{12}\,y_2 + \ldots + d_{1, p+1}\,y_{p+1}=\lambda y_1,\\
d_{21}\,y_1 + d_{22}\,y_2 + \ldots + d_{2,p+2}\,y_{p+2}=\lambda y_2,\\
\ldots\\
d_{p,1}\, y_1 + d_{p,2}\, y_2 + \ldots + d_{p,2p}\, y_{2p}=\lambda y_p,\\
d_{n,n-p}\, y_{n-p} + d_{n,n-p+1}\, y_{n-p+1} + \ldots + d_{n,
n+p}\, y_{n+p}=\lambda y_n, \quad n=p+1, p+2, \ldots.\\
\end{array}
\right.
\end{equation}
It is advisable to define coefficients $d_{i,j}$ for the indices
with $\min(i,j)\le 0$ as follows:
\begin{equation} \label{1.3}
d_{ij}=1, \quad |i-j|=p, \qquad d_{ij}=0, \quad |i-j|\not=p,
\end{equation}
and so system (\ref{1.2}) is equivalent to
\begin{equation} \label{1.4}
d_{n,n-p}\, y_{n-p} + d_{n, n-p+1}\, y_{n-p+1} + \ldots + d_{n,
n+p}\, y_{n+p}=\lambda y_n, \quad n\in\mathbb{N}
\end{equation}
with the initial conditions
\begin{equation} \label{1.5}
y_{1-p}=y_{2-p}=\ldots=y_0=0.
\end{equation}
Thus, the vector $\vec{y}=\{y_n\}_{n\ge 1-p}\in\ell^2$ is the
eigenvector of the operator $D$ corresponding to the eigenvalue
$\lambda$ if and only if $\{y_n\}$ satisfies (\ref{1.4}),
(\ref{1.5}).

It seems natural to analyze equation (\ref{1.4}) within the
framework of the general theory of linear difference equations.
The equation
\begin{equation} \label{1.6}
x(n+k) + a_1(n) x(n+k-1) + \ldots + a_k(n) y(n) = 0
\end{equation}
is said to belong to the Poincar\'e class if $a_k(n) \ne 0$ and
there exist limits (in $\mathbb{C}$)
$$ b_j = \lim_{n\to\infty} a_j(n), \quad j=1,2,\ldots,k.
$$
Denote by $\{w_j\}_{j=1}^k$ all the roots (counting the
multiplicity) of the characteristic equation
\begin{equation} \label{1.7}
w^k+b_1w^{k-1}+\ldots+b_k=0.
\end{equation}
One of the cornerstones of the theory of linear difference
equations is the following  result due to Perron.

\noindent {\bf Theorem} (\cite{Perron} , Satz 3). {\it Let the roots $\{w_j\}$
of $(\ref{1.7})$ lie on the circles $\Gamma_l=\{|w|=\rho_l\},\
l=1,2,\ldots,m,\ \rho_j\ne \rho_k$, and exactly $v_l\ge 1$ of them
$($counted according with their multiplicity$)$ belongs to each
circle $\Gamma_l$, so $\nu_1+\ldots+\nu_m=k$. Then $(2.6)$ has a
fundamental system of solutions
$$ S=\{y_1,\ldots,y_k\} = \bigcup_{l=1}^m S_l\, ,
$$
the sets $\{S_l\}$ are disjoint, $|S_l|=\nu_l$, and for any
nontrivial linear combination $y(n)$ of the solutions from $S_l$
\begin{equation} \label{1.8}
\limsup_{n\to\infty} \sqrt[n]{|y(n)|}=\rho_l,\qquad
l=1,2,\ldots,m,
\end{equation}
holds.}

\begin{proposition} \label{propos2.1} {\it For any $\lambda\in
\mathbb{C}\backslash[-2,2]$ the dimension of the space of
$\ell^2$-solutions of $(\ref{1.4})$ equals $p$.}
\end{proposition}

\noindent {\it Proof.} Note that equation (\ref{1.4}) has order
$k=2p$ (after dividing through by the leading coefficient $d_{n,\,
n+p}$) and belongs to the Poincar\'e class by assumption
(\ref{1.2}). Its characteristic equation (\ref{1.7}) has now the
form
$$ w^{2p} - \lambda w^p +1 = (w^p-z)(w^p-z^{-1})=0\, : \quad
\lambda=z+z^{-1}, \ {z}<1\,.
$$
For its roots we have
$$
|w_1|=\ldots=|w_p|=|z|<1<|z|^{-1}=|w_{p+1}|=\ldots=|w_{2p}|.
$$
By the Perron theorem, the fundamental system $S$ of solutions of
(\ref{1.4}) exists
$$ S=\{y_1,\ldots,y_p;y_{p+1},\ldots,y_{2p}\}=S_1 \cup S_2\,; \quad
\dim\, {\rm span} \, S_1 = \dim\, {\rm span} \, S_2 = p\,,
$$
and each solution $y\in {\rm span}\, S_1$ is in $\ell^2$ (and even
decreases exponentially fast). Let now $y$ be any solution of
(\ref{1.4}) from $\ell^2$,
$$ y=\sum_{j=1}^p c_j y_j + \sum_{j=p+1}^{2p} c_j y_j = y' + y''\,.
$$
But $y'\in\ell^2$, and so $y''\in\ell^2$ which by (\ref{1.8}) and
$|z|^{-1}>1$ is possible only when $c_j = 0$ for
$j=p+1,\ldots,2p$, as needed. \hfill $\square$

\begin{proposition} \label{propos2.2}
{\it Let $\{y_n^{(i)}\}_{n\ge 1-p}$, $ \
i=1,2,\ldots,p$, be linearly independent solutions of
$(\ref{1.4})$ from $\ell^2$. The number $\lambda$ is an eigenvalue
of the operator $D$ if and only if}
\begin{equation} \label{1.9}
\det Y_0(\lambda)=
\begin{vmatrix}
  y_{1-p}^{(1)} & y_{2-p}^{(1)} & \ldots & y_{0}^{(1)} \\
  \ldots & \ldots & \ldots & \ldots \\
  y_{1-p}^{(p)} & y_{2-p}^{(p)} & \ldots & y_{0}^{(p)}
\end{vmatrix}=0
\end{equation}
\end{proposition}

\noindent {\it Proof.} Suppose that $\det Y_0(\lambda)=0$. Then
there are numbers $\alpha^{(1)},\ldots,\alpha^{(p)}$ , which do
not vanish simultaneously, such that

$$\left\{
\begin{array}{ccc}
  \alpha^{(1)} y_{1-p}^{(1)} &+\ldots &+ \alpha^{(p)} y_{1-p}^{(p)} = 0\,,
\\
  \ldots & \ldots & \ldots \\
  \alpha^{(1)} y_0^{(1)} &+\ldots  &+ \alpha^{(p)} y_0^{(p)} = 0\,.
\end{array}
\right.
$$
Hence the linear combination
\begin{equation} \label{1.10}
y_n=\alpha^{(1)} y_n^{(1)} + \ldots + \alpha^{(p)} y_n^{(p)},
\qquad n\ge 1-p,
\end{equation}
belongs to $\ell^2$ and satisfies (\ref{1.4}), (\ref{1.5}), i.e.,
$\lambda$ is an eigenvalue of the operator $D$.

Conversely, let $\lambda$ be an eigenvalue and $y=\{y_n\}_{n\ge
1}$ the corresponding eigenvector. Then $\{y_n\}_{n\ge 1-p}$ is an
$\ell^2$-solution of (\ref{1.4}) with the initial conditions
(\ref{1.5}). By Proposition \ref{propos2.1} (\ref{1.10}) holds with
coefficients $\alpha^{(1)},\ldots,\alpha^{(p)}$ which do not
vanish simultaneously. Then (\ref{1.9}) follows immediately from
(\ref{1.5}). \hfill $\square$

\setcounter{section}{2} \setcounter{equation}{0}
\section {The matrix-valued Jost solution}

The goal of this Section is to establish the existence of
matrix-valued analogue of the Jost solution for the banded matrix
$D$. Once we have the matrix Jost solution at our disposal, we
will be able to construct $p$ linearly independent square-summable
solutions of (\ref{1.4}) and, in view of Proposition \ref{propos2.2}, to
reduce the study of the location of the discrete spectrum for the
matrix $D$ to the location of the zeros for determinant
(\ref{1.9}), composed of these $p$ solutions.

It is convenient to rewrite the initial equation in the form of a
three-term recurrence matrix relation, by looking at $D$ as a
block-Jacobi matrix. Along this way we can extend the standard
techniques of proving the existence of the Jost solution for
Jacobi matrices to the case of banded matrices.

Define the following $p\times p$-matrices:

\begin{equation} \label{2.1}
\begin{aligned} 
  A_k &= \left(\begin{array}{ccc} d_{(k-1)p+1,\, (k-2)p+1} &  \ldots &
d_{(k-1)p+1,\, (k-1)p} \\ \vdots &  & \vdots\\
d_{kp,\, (k-2)p+1} & \ldots & d_{kp,\, (k-1)p}\\
\end{array}\right), \\
  B_k &= \left(\begin{array}{ccc} d_{(k-1)p+1,\, (k-1)p+1} & \ldots &
d_{(k-1)p+1,\, kp} \\ \vdots & \ & \vdots\\
d_{kp,\, (k-1)p+1} & \ldots & d_{kp,\, kp}\\
\end{array}\right), \\
  C_k &= \left(\begin{array}{ccc} d_{(k-1)p+1,\, kp+1} & \ldots &
d_{(k-1)p+1,\, (k+1)p} \\ \vdots & \ & \vdots\\
d_{kp,\, kp+1} & \ldots & d_{kp,\, (k+1)p}\\
\end{array}\right).
\end{aligned}
\end{equation}

Then the matrix $D$ can be represented in the form
\begin{equation} \label{2.2}
D = \left(\begin{array}{ccccc}
B_1 & C_1 & 0 & 0 & \ldots \\
A_2 & B_2 & C_2 & 0 &
\ldots \\
0 & A_3 & B_3 & C_3 & \ldots \\
0 & 0 & A_4 & B_4 & \ldots\\
\ldots & \ldots & \ldots & \ddots & \ddots\\
\end{array}\right),
\end{equation}
with the upper triangular matrices $A_k$ and the lower triangular
$C_k$ which are invertible due to (\ref{0.1}). To be consistent
with (\ref{1.3}) we put $A_1=C_0=I$ a unit $p\times p$ matrix,
$B_0=0$.

Having $p$ solutions $\{\varphi_j^{(l)}\}_{j\ge 1-p}$,
$l=1,2,\ldots,p$ of equation (\ref{1.4}) at hand, we can make up
$p\times p$ matrices
\begin{equation} \label{2.3}
\Psi_j = \left(\begin{array}{ccc} \varphi_{(j-1)p+1}^{(1)} &
\ldots &
\varphi_{(j-1)p+1}^{(p)} \\ \vdots &  & \vdots\\
\varphi_{jp}^{(1)} & \ldots & \varphi_{jp}^{(p)}\\
\end{array}\right)\,, \quad j\in\mathbb{Z}_+:=0,1,\ldots.
\end{equation}
and write (\ref{1.4}) in the matrix form
\begin{equation} \label{2.4}
A_k\Psi_{k-1} + B_k\Psi_k + C_k\Psi_{k+1} = \lambda \Psi_k, \quad
k\in\mathbb N.
\end{equation}

It will be convenient to modify equation (\ref{2.4}), getting rid
of the coefficients $A_k$'s. Suppose that
\begin{equation} \label{2.5}
\sum_{k=1}^\infty \| I-A_k \| <\infty \,,
\end{equation}
where $\|\cdot\|$ is any norm in the space of matrices. It is well
known, that there exists an infinite product (from the right to
the left)
$$ A:= \prod_{j=1}^\infty A_j =
\lim_{n\rightarrow\infty} (A_k A_{k-1} \ldots A_1) \,,
$$
and all the matrices $A$ and $A_n$ are invertible. Denote
\begin{equation} \label{2.5l}
L_j:=\prod_{i=j+1}^\infty A_j, \qquad L_jA_j = L_{j-1}, \qquad
\lim_{j\to\infty}L_j=I. \end{equation} The multiplication of
(\ref{2.4}) from the left by $L_k$ gives
$$ L_{k-1} \Psi_{k-1} + L_k B_k \Psi_k + L_k C_k \Psi_{k+1} =
\lambda L_k \Psi_k \,,
$$
$$ L_{k-1} \Psi_{k-1} + L_k B_k L_k^{-1} \cdot L_k \Psi_k + L_k C_k
L_{k+1}^{-1} L_{k+1} \Psi_{k+1}= \lambda L_k \Psi_k \,.
$$
Hence the matrices
\begin{equation}\label{2.52}
\Phi_k:=L_k\Psi_k
\end{equation}
satisfy
\begin{equation} \label{2.6}
\Phi_{k-1} + \tilde B_k \Phi_k + \tilde C_k \Phi_{k+1} = \lambda
\Phi_k \,, \qquad k\in\mathbb{N}
\end{equation}
with
\begin{equation} \label{2.6'}
\tilde{B_k} = L_k B_k L_k^{-1},  \qquad \tilde{C_k}:= L_k C_k
L_{k+1}^{-1} \,.
\end{equation}

For the definiteness sake we choose the ``row norm''
$$ \|T\|:=\max_{1\le k\le p} \sum_{j=1}^p |t_{kj}|, \qquad
T=\{t_{kj}\}_{k,j=1}^p. $$ Then by (\ref{2.1})
\begin{equation} \label{2.5x}
\max\bigl(\|A_k-I\|, \|B_k\|, \|C_k-I\|\bigr)\leq \hat
q_k:=\max_{1\le j\le p} (q_{(k-1)p+j}), \quad k\in\mathbb N,
\end{equation}
$q_k$ are defined in (\ref{0.3}). In accordance with (\ref{1.3})
$q_{1-p}=\ldots=q_0=0$, so we put $\hat q_0=0$. It is clear that
\begin{equation} \label{2.5y}
\sum_{k=1}^\infty \hat q_k\le \sum_{k=1}^\infty q_k\le
p\sum_{k=1}^\infty \hat q_k
\end{equation}
so (\ref{2.5}) holds whenever $\{q_n\}\in\ell^1$.

We will use the complex parameter $z$ related to the spectral
parameter $\lambda$ by the Zhukovsky transform:
$$ \lambda=z+z^{-1}; \quad |z|<1\,.
$$

Denote by $g$ the Green kernel

\begin{equation} \label{2.8}
g(n,k,z)=\left\{ \begin{array}{cc} {\displaystyle
\frac{z^{k-n}-z^{n-k}}{z-z^{-1}} },&k>n,\\0,&k\leq n,
\end{array}\right.\quad n,k\in\mathbb{Z}_+:=\{0,1,...\},\ z\neq 0.
\end{equation}
It is clear that $g(n,k,z)$ satisfies the recurrence relations
\begin{eqnarray} \label{2.9}
g(n,k+1,z)+g(n,k-1,z)-(z+z^{-1}) g(n,k,z)&= \delta(n,k), \\
\label{2.10} g(n-1,k,z)+g(n+1,k,z)-(z+z^{-1}) g(n,k,z)&=
\delta(n,k),
\end{eqnarray}
where $\delta(n,k)$  is the Kronecker symbol.

We proceed with the following conditional result.

\begin{proposition} Suppose that equation $(\ref{2.6})$ has a solution
$V_n$ with the asymptotic behavior at infinity
\begin{equation}\label{2.11}
\lim_{n\to\infty}\,V_n(z) z^{-n} = I
\end{equation}
for $z\in\mathbb{D}$. Then $V_n$ satisfies the discrete integral
equation
\begin{equation} \label{2.12}
V_n(z)=z^n I + \sum_{k=n+1}^\infty J(n,k,z)\, V_k(z),\quad
n\in\mathbb{N},
\end{equation}
with
\begin{equation}\label{2.13}
J(n,k,z)=- g(n,k,z) \tilde{B}_k + g(n,k-1,z) \left(I - \tilde
C_{k-1}\right).
\end{equation}
\end{proposition}

{\it Proof.} Let us multiply (\ref{2.9}) by $V_k$, (\ref{2.6}) for
$V_k$ by $g(n,k)$, and subtract the latter from the former:
$$
 g(n,k+1)V_k + g(n,k-1) V_k - g(n,k)V_{k-1} - g(n,k)\tilde B_k
 V_k - g(n,k)\tilde C_k V_{k+1} = \delta(n,k) V_k.
$$
Summing up over $k$ from $n$ to $N$ gives
$$
\begin{aligned} V_n &= \sum_{k=n+1}^N\,\left\{-g(n,k)\tilde B_k +
g(n,k-1)\left(I - \tilde C_{k-1}\right) \right\}V_k \\
&+ g(n,N+1)V_N - g(n,N+1)\tilde C_N V_{N+1}.
\end{aligned}
$$ For $|z|<1$ we
have by (\ref{2.8}) and (\ref{2.11})
$$\lim_{N\to\infty}\,\bigl( g(n,N+1)V_N - g(n,N)\tilde C_N V_{N+1}\bigl) =
z^n I,$$ which along with $J(n,n)=0$ leads to (\ref{2.12}), as
needed.\hfill $\square$

The converse statement is equally simple.

\begin{proposition} Each solution $\{V_n(z)\}_{n\geq0}$,
$z\in\overline{\mathbb{D}}$, of equation $(\ref{2.12})$ with
$n\in\mathbb{Z}_+$ satisfies the three-term recurrence relation
$(\ref{2.6})$.
\end{proposition}

{\it Proof.} Write for $n\geq 1$
$$
\begin{aligned}V_{n-1} + V_{n+1} &= (z^{n-1} + z^{n+1}) I +
\sum_{k=n}^\infty\,J(n-1,k)V_k +\sum_{k=n+2}^\infty\,J(n+1,k)V_k
\\ &= (z+z^{-1})z^n I + J(n-1,n) V_n + J(n-1,n+1)V_{n+1} \\
&+ \sum_{k=n+2}^\infty\,\{J(n-1,k) +J(n+1,k)\}V_k. \end{aligned}
$$
By (\ref{2.8}), (\ref{2.13}) and (\ref{2.10})
$$J(n-1,n)=-\tilde B_n, \qquad J(n-1,n+1)=-(z+z^{-1})\tilde B_{n+1}+I-\tilde
C_n
$$
and
$$ J(n-1,k)+J(n+1,k)=(z+z^{-1})J(n,k), \qquad k\ge n+2.
$$
Hence
$$ \begin{aligned} &{} V_{n-1}+V_{n+1} + \tilde B_n V_n -(I - \tilde C_n)
V_{n+1} = (z+z^{-1}) z^n I - \tilde B_n V_n - (z+z^{-1}) \tilde
B_n V_{n+1} \\ &+ (I-\tilde C_n) V_{n+1} + \sum_{k=n+2}^\infty
\left\{ (z+z^{-1} J(n,k)\right\} V_k +
\tilde B_n V_n - (I-\tilde C_n) V_{n+1} \\
&= (z+z^{-1}) \left(z^n +\sum_{k=n+1}^\infty\,J(n,k)v_k\right) =
(z+z^{-1})V_n, \end{aligned}$$
 which is exactly (\ref{2.6}).\hfill $\square$

{\bf The Jost solution.} To analyze equation (\ref{2.12}) we
introduce new variables
\begin{equation} \label{2.13.1}
\tilde V_n(z):= z^{-n} V_n, \quad \tilde J(n,k,z):=z^{k-n}
J(n,k,z),
\end{equation}
so that, instead of (\ref{2.12}), we have
\begin{equation}\label{2.14}
\tilde V_n(z) = I + \sum_{k=n+1}^\infty\,\tilde J(n,k,z) \tilde
V_k(z), \quad n\in\mathbb{Z}_+.
\end{equation}
Now $\tilde J(n,m,\cdot)$ is a polynomial with matrix
coefficients. Since
$$|g(n,k,z)z^{k-n}|=\frac{|z^{2(k-n)}-1|}{|z-z^{-1}|}
\leq |z|\min\left\{|k-n|,\frac2{|z^2-1|}\right\},
$$
the kernel $\tilde J$ is bounded by
\begin{equation}\label{2.15}
\|\tilde J(n,k,z) \| \leq |z|
\min\left\{|k-n|,\frac2{|z^2-1|}\right\}\,h_k,\quad
z\in\overline{\mathbb{D}}
\end{equation}
where
\begin{equation}\label{2.16}
h_k:=\|\tilde B_k \| + \| I- \tilde
C_{k-1}\|=\|L_kB_kL_k^{-1}\|+\|I-L_{k-1}C_{k-1}L_k^{-1}\|, \quad
k\in\mathbb N
\end{equation}
(see (\ref{2.6'})). We have
\begin{equation} \label{2.16x}
h_k\le
\|L_k\|\|L_k^{-1}\|\left(\|B_k\|+\|I-A_k\|+\|A_k\|\|I-C_{k-1}\|\right).
\end{equation}
Since $\|A_k\|\le C(D)$ and by (\ref{2.5l})
$\|L_k\|\cdot\|L_k^{-1}\|\le C(D)$ (throughout the rest of the
paper $C=C(D)$ stands for various positive constants which depend
only on $p$ and the original matrix $D$), we see from (\ref{2.5x})
that
\begin{equation} \label{2.16.1}
h_k\le C(D)(\hat q_{k-1}+\hat q _k).
\end{equation}

The existence of the Jost solutions for equation (\ref{2.6}) will
be proved under the assumption
$$ \sum_i q_i=\sum_i \left(|d_{i,\, i-p} - 1| + | d_{i,\, i+p}-1 |
+\sum_{r=1-p}^{p-1} |d_{i,\, i+r}|\right) <\infty \,,
$$
which by (\ref{2.5x}), (\ref{2.5y}) and (\ref{2.16.1}) implies
\begin{equation} \label{2.7}
\sum_k (\|I-A_k\| + \|B_k\| + \|I-C_k\| ) <\infty \,, \qquad
\sum_k h_k < \infty.
\end{equation}

The main result concerning equation (\ref{2.12}) is the following

\begin{theorem} $(i)$ Suppose that
\begin{equation}\label{2.17}
\sum_{k=1}^\infty q_k\,<\,\infty \,,
\end{equation}
$q_k$ are defined in $(\ref{0.3})$. Then equation $(\ref{2.12})$
has a unique solution $V_n$, which is analytic in $\mathbb{D}$,
continuous on
$\mathbb{D}_1:=\overline{\mathbb{D}}\setminus\{\pm1\}$ and
\footnote{Following the terminology of self-adjoint case for
Jacobi matrices, we call this solution the {\it Jost solution}.
The function $V_0$ is the matrix analogue of the {\it Jost
function}.}
\begin{equation}\label{2.18}
\|V_n - z^n I\|\leq C|z|^n \left\{\frac{|z|}{|z^2 - 1|}
\sum_{k=n}^\infty\,
 q_k\right\}\,\exp\left\{\frac{C|z|}{|z^2-1|}
\sum_{k=n}^\infty\,q_k\right\}
\end{equation}
for $z\in\mathbb{D}_1$, $n\in\mathbb{Z}_+$.

$(ii)$ Suppose that
\begin{equation}\label{2.19}
\sum_{k=1}^\infty k q_k\,<\,\infty.
\end{equation}
Then $V_n$ is analytic in $\mathbb{D}$, continuous on
$\overline{\mathbb{D}}$ and
\begin{equation}\label{2.20}
\|V_n - z^n I \|\leq C|z|^n\left\{\sum_{k=n}^\infty\, k
q_k\right\}\, \exp\left\{C\sum_{k=n}^\infty\,k q_k\right\},\ \
z\in\mathbb{D},\quad n\in\mathbb{Z}_+.
\end{equation}
\end{theorem}
{\it Proof.} We apply the method of successive approximations.
Write (\ref{2.14}) as
\begin{equation} \label{2.21}
F_n(z) = G_n(z) + \sum_{k=n+1}^\infty\,\tilde J(n,k,z) F_k(z)
\end{equation}
with
\begin{equation}\label{2.22}
F_k(z):=\tilde V_k (z) - 1 , \quad G_n(z):=\sum_{k=n+1}^\infty
\,\tilde J(n,k,z).
\end{equation}
\noindent (i) By (\ref{2.15})
\begin{equation}\label{2.23}
\|\tilde J(n,k,z)\|\leq\phi(z)\,h_k,\quad z\in\mathbb{D}_1, \qquad
\phi(z):=2|z||z^2-1|^{-1}.
\end{equation}
The series in (\ref{2.22}) converges uniformly on compact subsets
of $\mathbb{D}_1$ by (\ref{2.17}), (\ref{2.7}), and so $G_n$ is
analytic in $\mathbb{D}$ and continuous on $\mathbb{D}_1$.  As a
starting point for the method of successive approximation, we put
$F_{n,1} = G_n$ and denote
$$ F_{n,j+1}(z) := \sum_{k=n+1}^\infty\,\tilde J(n,k,z)F_{k,j}(z).
$$
Let $\sigma_0(n):=\sum_{k=n+1}^\infty\,h_k$. By induction on $j$
we prove that
\begin{equation}\label{2.24}
\left\| F_{n,j}(z)\right\| \leq\,
\frac{(\phi(z)\sigma_0(n))^j}{(j-1)!}.
\end{equation}
Indeed, for $j=1$ we have $F_{n,1}=G_n$ and the result holds by
the definition of $\sigma_0$ and (\ref{2.23}). Next, let
(\ref{2.24}) be true. Then
$$|F_{n,j+1}(z)| \leq \phi(z)\,\sum_{k=n+1}^\infty\,h_k \|F_{k,j}(z)\| \leq
\frac{(\phi(z))^{j+1}}{(j-1)!}\sum_{k=n+1}^\infty\, h_k
\sigma_0^j(k).
$$
An elementary inequality $(a+b)^{j+1} - a^{j+1}\geq(j+1)ba^j$
gives
$$\sum_{k=n+1}^\infty\,h_k\sigma_0^j(k)\leq\frac{1}{j}
\sum_{k=n+1}^\infty\{\sigma_0^{j+1}(k-1)
-\sigma_0^{j+1}(k)\}=\frac{\sigma_0^{j+1}(n)}{j}\,,
$$
which proves (\ref{2.24}) for $F_{n,j+1}$. Thereby the series
$$F_n(z)=\sum_{j=1}^\infty\,F_{n,j}(z)$$
converges uniformly on compact subsets of $\mathbb{D}_1$ and
solves (\ref{2.21}), being analytic in $\mathbb{D}$ and continuous
on $\mathbb{D}_1$. It is also clear from (\ref{2.24}) that
\begin{equation} \label{2.24x}
\|F_n(z)\|=\|\tilde V_n(z)-I\|\le\sum_{j=1}^\infty
\|F_{n,j}(z)\|\le \phi(z)\sigma_0(n)\exp\{\phi(z)\sigma_0(n)\}.
\end{equation}
To reach (\ref{2.18}) it remains only to note that by
(\ref{2.16.1})
\[
\sigma_0(n)\le C\sum_{k=n+1}^\infty (\hat q_{k-1}+\hat q_k)\le
2C\sum_{k=n}^\infty \hat q_k\le 2C\sum_{j=n}^\infty q_j
\]
(the latter inequality easily follows from the definition of $\hat
q_k$).

To prove uniqueness suppose that there are two solutions $F_n$ and
$\tilde F_n$ of (\ref{2.21}). Take the difference and apply
(\ref{2.23}):
\begin{equation}
\begin{array}{c}
{\displaystyle | F_n - \tilde F_n (z) | = | \sum_{k=n+1}^\infty
\tilde J(n,k,z) \left( F_k(z)-\tilde F_k(z) \right) | } \\
{\displaystyle s_n\leq \sum_{m=n+1}^\infty \phi(z) s_m h_m=r_n, }\\
\end{array} \label{2.25}
\end{equation}
where $s_n:=\|F_n(z) -\tilde F_n(z)\|$.

Clearly, $r_n\to 0$ as $n\to\infty$ and if $r_m=0$ for some $m$,
then by (\ref{2.25}) we have $s_n\equiv 0$. If $r_n>0$, then
\begin{equation} \label{2.26}
\frac{r_{n-1}-r_n}{r_n}=\frac{s_n\phi(z)h_n}{r_n}\leq \phi(z)h_n,
\quad r_k\leq\prod_{j=k+1}^M\left(1+\phi(z)h_j\right)\,r_M
\end{equation}
which leads to $r_m=0$ and again $s_n\equiv 0$. So the uniqueness
is proved.

\noindent (ii) The same sort of reasoning is applicable with
$$\|\tilde J(n,k,z) \| \leq |z| |k-n| h_k \leq k h_k
$$
and
\begin{equation} \label{2.26z}
\| F_{n,j}(z) \| \leq \frac{\sigma_1^j(n)}{(j-1)!}\,,\quad
\sigma_1(n):=\sum_{k=n+1}^\infty k h_k
\end{equation}
instead of (\ref{2.23}) and (\ref{2.24}), respectively. We have
\begin{equation} \label{2.26a}
\|F_n(z)\|=\|\tilde V_n(z)-I\|\le\sum_{j=1}^\infty
\|F_{n,j}(z)\|\le \sigma_1(n)\exp\{\sigma_1(n)\}
\end{equation}
and
\[
\sigma_1(n)\le C\sum_{k=n+1}^\infty k(\hat q_{k-1}+\hat q_k)\le
2C\sum_{k=n}^\infty k\hat q_k\le 2C\sum_{j=n}^\infty j q_j
\]
(the latter inequality easily follows from the definition of $\hat
q_k$).                      \hfill$\square$

\medskip

{\sl Remark}. The constants $C$ that enter (\ref{2.16.1}),
(\ref{2.18}) and (\ref{2.20}) are inefficient. This circumstance
makes no problem when studying the limit set for the discrete
spectrum. In contrast to this case, the efficient constants are
called for when dealing with domains which contain the whole
discrete spectrum. Such constants will be obtained in the next
section under additional assumptions of ``non asymptotic flavor''.

Throughout the rest of the section we assume that condition
(\ref{2.19}) is satisfied. It is clear that equation (\ref{2.6})
can be rewritten for the functions $\tilde V_n(z)$, defined in
(\ref{2.13.1}), as
\begin{equation}\label{2.27}
\begin{aligned}
\tilde V_n(z) &= (\lambda - \tilde B_n) z \tilde V_{n+1}(z) -
\tilde C_n z^2 \tilde V_{n+2}(z) \\ &= (z^2+1-\tilde B_n z) \tilde
V_{n+1}(z) - \tilde C_n z^2 \tilde V_{n+2}(z) \,.
\end{aligned}
\end{equation}

Let us now expand $\tilde V_n(z)$ in the Taylor series taking into
account definition (\ref{2.13.1}) and (\ref{2.18})
\begin{equation}\label{2.28}
\tilde V_n(z) = I + \sum_{j=1}^\infty K(n,j)z^j \,.
\end{equation}
Here $\|K(n,j)\|_{n,j=1}^\infty$ is the operator which transforms the Jost solutions
of (\ref{1.1}) for $D=D_0$ ti that of (\ref{1.1}) for $D$.
If we plug (\ref{2.28}) into (\ref{2.27}) and match the
coefficients for the same powers $z^j$ we have
$$ j=1: \quad K(n,1) = K(n+1,1) - \tilde B_n \,,
$$
$$ j=2: \quad K(n,2) = I + K(n+1,2) - \tilde B_n K(n+1,1) - \tilde
C_n \,,
$$
$$ j \ge 2: \quad K(n,j+1)=K(n+1,j-1) + K(n+1,j+1) - \tilde B_n
K(n+1,j) - \tilde C_n K(n+2,j-1).
$$
Summing up each one of these expressions for $k=n,n+1,\ldots$, it
is not hard to verify that
\begin{equation}\label{2.29}
K(n,1)= - \sum_{k=n+1}^\infty \tilde B_{k-1} \,,
\end{equation}
\begin{equation}\label{2.30}
K(n,2)=-\sum_{k=n+1}^\infty \left\{\tilde B_{k-1} K(k,1) +(\tilde
C_{k-1} - I)\right\}\,, \end{equation}
\begin{equation}\label{2.31}
K(n,j+1)=K(n+1,j-1) - \sum_{k=n+1}^\infty \left \{\tilde B_{k-1}
K(k,j)+\left(\tilde C_{k-1} - I\right) K(k+1,j-1)\right\}
\end{equation}
In the last step we used $K(n,j)\to 0$ for $n\to\infty$ and any
fixed $j$, which follows from the Cauchy inequality and
(\ref{2.20})
$$ \| K(n,j) \| \leq \max \| \tilde V_n(z) - I\|
\leq C \sum_{k=n}^\infty kq_k \,.
$$

From (\ref{2.29})--(\ref{2.31}), using the induction on $j$, we
obtain
\begin{equation} \label{2.32}
\| K(n,j) \| \leq \kappa(n,j)
\kappa\left(n+\left[\frac{j}{2}\right]\right), \quad
n\in\mathbb{Z}_+,
\end{equation}
where $\kappa(n)$ and $\kappa(n,m)$ are defined by
\begin{equation} \label{2.33}
\kappa(n):=\sum_{j=n}^\infty g_j, \quad
\kappa(n,m):=\prod_{j=n+1}^{n+m-1}(1+\kappa(j))
=\prod_{j=1}^{m-1}(1+\kappa(n+j)) \,,
\end{equation}
and
\begin{equation} \label{2.333}
g_j = \|\tilde B_j\| + \|I - \tilde C_j \| .
\end{equation}

\medskip

In fact, for $j=1$ we have $\kappa(n,1)=1$, $\left[\frac{1}{2}\right]=0$ and $\|K(n,1)\| \leq
\kappa(n+1) \leq \kappa(n)$. Further, for $j=2$,
$$ \begin{array}{ll}
\|K(n,2) \| & \leq {\displaystyle \sum_{k=n+1}^\infty \left\{\|\tilde B_{k-1}\| \|K(k,1)\|
+\|\tilde C_{k-1} - I\| \right\}  }\\
& \leq {\displaystyle \sum_{k=n+1}^\infty \left\{ \left( \|K(k,1)\|+1\right)
\left(  \|\tilde B_{k-1}\| + \|\tilde C_{k-1} - I\| \right) \right\| }\\
& \leq {\displaystyle \sum_{k=n+1}^\infty g_k \left( 1+ \kappa(k) \right) } \\
& \leq {\displaystyle \left( 1+ \kappa(n+1)\right) \kappa(n+1) } \\
& = {\displaystyle \kappa(n,2) \kappa(n+1). } \\
\end{array}
$$
When we pass from the even $j=2l$ to the odd $2l+1$, we have by
(\ref{2.31}) and by the inductive hypothesis:
$$ \| K(n, 2l+1) \| \leq \kappa(n+1, 2l-1)) \kappa(n+l) +
\sum_{k=n+1}^\infty \left \{\|\tilde B_{k-1}\|
\|K(k,2l)\|+ \|\tilde C_{k-1} - I\| \| K(k+1,2l-1)\| \right\} .
$$
But according to the inductive hypothesis
$$ K(k,2l)\| \leq \kappa(k,2l)\kappa(k+l)\,;
$$
$$ \| K(k+1, 2l-1) \| \leq \kappa(k+1, 2l-1) \kappa(k+l) \leq \kappa(k,2l)\kappa(k+l).
$$
Using these inequalities, obtain
$$ \begin{array}{ll}
\|K(n,2l+1) \| & \leq {\displaystyle \kappa(n+1,2l-1)\kappa(n+l)
+ \sum_{k=n+1}^\infty d_k \kappa(k,2l) \kappa(k+l) }\\
& \leq \kappa(n+1,2l-1)\kappa(n+l)
+ \kappa(n+l) \kappa(n+1,2l) \kappa(n+1) \\
& = \kappa(n+l) \left\{ \kappa(n+1,2l-1)
+ \kappa(n+1,2l) \kappa(n+1) \right\} \\
\end{array}
$$
But
$$ \begin{array}{l}
\kappa(n+1,2l-1) + \kappa(n+1,2l) \kappa(n+1) = \\
{\displaystyle = \prod_{j=n+2}^{j=n+2l-1} \left(1+\kappa(j)\right) +
\prod_{j=n+2}^{j=n+2l} \left(1+\kappa(j)\right) \kappa(n+1) } \\
{\displaystyle = \prod_{j=n+2}^{j=n+2l-1} \left(1+\kappa(j)\right)
\left\{ 1+ \left(1+\kappa(n+2l)\right) \kappa(n+1) \right\} }\\
{\displaystyle \leq \prod_{j=n+2}^{j=n+2l-1} \left(1+\kappa(j)\right)
\left\{ 1+ \kappa(n+2l) + \left(1+\kappa(n+2l)\right) \kappa(n+1) \right\} }\\
{\displaystyle = \prod_{j=n+2}^{j=n+2l-1} \left(1+\kappa(j)\right)
\left\{ \left( 1+ \kappa(n+2l)\right) \left(1+\kappa(n+2l)\right) \right\} } \\
= \kappa(n,2l+1),\\
\end{array}
$$
from which we have
$$ \|K(n,2l+1)\| \leq \kappa(n+2l) \kappa(n+l),
$$
as needed. Analogous calculations help us to pass from the odd $j=2l+1$ to
the even $j+1=2l+2$.

\medskip

It is easy to see that $\kappa(n)$ and $\kappa(n,m)$ in (\ref{2.33}) can be replaced by
\begin{equation} \label{2.334}
\tilde \kappa(n):= C \sum_{j=n}^\infty q_j, \quad
\tilde \kappa(n,m):=\prod_{j=n+1}^{n+m-1}(1+\tilde\kappa(j))
=\prod_{j=1}^{m-1}(1+\tilde\kappa(n+j)) \,.
\end{equation}
Further, it is evident that $\{\kappa(n)\} \in \ell^1$ and the
sequences $\tilde \kappa(\cdot)$ and $\tilde \kappa(\cdot,m)$ decrease
monotonically. Hence
\begin{equation}
\|K(n,j)\| \leq \prod_{j=1}^{\infty}(1+\tilde \kappa(j)) \,
\sum_{k=n+\left[\frac{j}{2}\right]}^\infty \, q_k.
\end{equation}
Taking the latter expression with $n=0$, we come to the following

\begin{theorem} \label{theor3.4}
Under hypothesis $(\ref{2.19})$ the
Taylor coefficients of the matrix-valued function
\begin{equation}\label{2.34}
\Delta(z):=V_0(z)=\sum_{j=0}^\infty \delta(j) z^j
\end{equation}
admit the bound
\begin{equation}\label{2.35}
\|\delta(j)\| \leq C\prod_{j=1}^{\infty}(1+\tilde\kappa(j))
\,\sum_{k=\left[\frac{j}{2}\right]}^\infty q_k,
\end{equation}
where
$[x]$ is an integer part of $x$. In particular, $\Delta$ belongs
to the space $W_+$ of absolutely convergent Taylor matrix-valued
series.
\end{theorem}

Let us now denote by $\Delta_{ij}$ the entries of the matrix
$\Delta$:
$$ \Delta := \|\Delta_{ij}\|_{i,j=1}^p \,.
$$

\begin{corollary} \label{coroll3.5} Let for the banded matrix $D$
the numbers
\begin{equation}\label{2.36}
M_{n+1}:=\sum_{k=0}^{\infty} (k+1)^{n+1} q_k < \infty.
\end{equation}
Then  the $n$'th derivative $\Delta^{(n)}(z)=V_0^{(n)}(z)$ belongs
to $W_+$ and
\begin{equation} \label{2.37}
\max_{z\in\overline{\mathbb{D}}} \left|\Delta^{(n)}_{ij}(z)\right|
\leq C(D)\,\frac{4^n}{n+1}\,M_{n+1}, \quad i,j=1,2\ldots,p.
\end{equation}
\end{corollary}

\noindent {\it Proof.} The statement is a simple consequence of
(\ref{2.36}), the series expansion
$$ \Delta^{(n)}(z)=\sum_{j=0}^\infty (j+1)\ldots
(j+n) \delta(j+n) z^j
$$
bounds (\ref{2.35}) and the obvious inequality
$|\Delta^{(n)}_{ij}(z)| \leq \| \Delta^{(n)}\|$. \hfill $\square$

\setcounter{section}{3} \setcounter{equation}{0}
\section {The limit set and location of the discrete spectrum}

Consider the matrix-valued Jost solutions $V_k(z)$, which exist
$\forall z \in\mathbb{C}\backslash\{0\}$. The $p$ \ scalar
solutions $\{v_j^{(l)}\}_{j\ge1-p}$, $l=1,2,\ldots,p$ of equation
(\ref{1.4}), constructed from $V_k(z)$ by formulae (\ref{2.52} and (\ref{2.3}), are
linearly independent and belong to $\ell^2$ due to asymptotic
formulae (\ref{2.18}). According to Proposition \ref{propos2.2}, the number $
\eta=\zeta +  \zeta^{-1}$ is an eigenvalue for the operator $D$ if and
only if the determinant of the matrix-valued function
$\Delta(z)=V_0(z)$ vanishes at the point $\zeta$. Thus, the study of
the discrete spectrum of the operator $D$ is reduced to the study
of the zeros of the function
\begin{equation}\label{3.1}
\gamma(z) := \det \Delta(z) .
\end{equation}

The main item of business in this section is the limit set $E_D$
of the discrete spectrum of the operator $D$. Remind that $E_D
\subset [-2,2]$.

Let $D\in\mathcal{P}_p(\beta)$ (see Definition \ref{def1}). Since
$\| L_k \|$ and $\|L_k^{-1}\|$ are uniformly bounded, it is clear
from $\ref{2.7}$ that $h_n$, defined in (\ref{2.16}), satisfies
the same inequality:
\begin{equation}\label{3.8x}
h_n \leq C_1\,\exp(-C_2\, (n+1)^\beta)
\end{equation}
(with the same exponent $\beta$, but other constants $C_1, C_2
>0$).

We are in a position now to prove the first result announced in
the introduction.

\begin{theorem} Let
$D\in\mathcal{P}(\beta)$ where $0<\beta<\frac{1}{2}$. Then $E_D$
is a closed point set of the Lebesgue measure zero and its
convergence exponent satisfies
\begin{equation} \label{3.8z} \dim E_D \leq
\tau(E_D)\leq\frac{1-2\beta}{1-\beta}\,,
\end{equation} where $\dim E_D$ is the Hausdorff dimension of
$E_D$. Moreover, if $J\in \mathcal{P}(\frac{1}{2})$ then
$E_D=\emptyset$, i.e., the discrete spectrum is finite.
\end{theorem}
{\it Proof.} \ Denote by $\mathcal{A}$ the set of all functions,
analytic inside $\mathbb D$ and continuous in $\overline{\mathbb
{D}}$. Recall that the set $E$ on the unit circle $\mathbb{T}$ is
called a zero set for a class $\mathcal{X}\subset\mathcal{A}$ of
functions, if there exists a non-trivial function
$f\in\mathcal{X}$, which vanishes on $E$. We want to show that the
function $\gamma$ belongs to a certain class $\mathcal{X}$ (the
Gevr\'e class, see below) with known properties of its zero sets.
Note that, since $\mathcal{X}\subset\mathcal{A}$ then, according
to the Fatou theorem, the zero set has the Lebesgue measure zero.

We begin with certain bounds for the derivatives of the function
$\gamma$, which can be obtained from Theorem \ref{theor3.4} and Corollary
\ref{coroll3.5}. To this end we write
$$ \gamma(z) := \det \Delta(z) = \sum_\pi {\rm sign} \, \pi \,
\Delta_{1,\pi(1)} \Delta_{2,\pi(2)} \ldots \Delta_{p,\pi(p)} \,,
$$
where $\pi$ are the permutations of the set $\{1,2,\ldots,p\}$.
For the $n$-th derivative of~$\gamma$
\begin{equation}\label{3.8}
\gamma^{(n)}(z) = \sum_\pi{\rm sign}\, \pi
\sum_{\substack{\sum_{1}^p k_j=n
\\
k_j\ge 0}}
\left(\begin{array}{c} n \\ k_1, k_2, \ldots , k_p \\
\end{array}\right) \Delta_{1,\pi(1)}^{(k_1)}
\Delta_{2,\pi(2)}^{(k_2)} \ldots \Delta_{p,\pi(p)}^{(k_p)}
\end{equation}
holds, where the multinomial coefficients
$\left(\begin{array}{c} n \\ k_1, k_2, \ldots , k_p \\
\end{array}\right)$ are defined from the identity
\begin{equation} \label{3.9}
\left(x_1 + x_2 + \ldots + x_p \right)^p =
\sum_{\substack{\sum_{1}^p k_j=n \\
k_j\ge 0}}
\left(\begin{array}{c} n \\ k_1, k_2, \ldots , k_p \\
\end{array}\right) x_1^{k_1} x_2^{k_2} \ldots x_p^{k_p}.
\end{equation}
(\ref{3.9}) with $x_1=x_2=\ldots=x_p=1$ gives
$$
\sum_{\substack{\sum_{1}^p k_j=n \\
k_j\ge 0}}
\left(\begin{array}{c} n \\ k_1, k_2, \ldots , k_p \\
\end{array}\right) = p^n \,.
$$

From (\ref{3.8}) and (\ref{2.37}) we immediately derive
\begin{equation} \label{3.10}
\max_{z\in\overline{\mathbb{D}}} \left|\gamma^{(n)}(z)\right|\leq
C^p p! \, 4^n \sum_{\substack{\sum_{1}^p k_j=n \\
k_j\ge 0}}
\left(\begin{array}{c} n \\ k_1, k_2, \ldots , k_p \\
\end{array}\right)\,M_{k_1+1}M_{k_2+1}\ldots M_{k_p+1},
\end{equation}
where the numbers $M_r$ are defined in (\ref{2.36}). To estimate
the product $M_{k_1+1} \ldots M_{k_p+1}$ note first, that by
(\ref{3.8x})
$$ M_r\leq C\sum_{k=0}^\infty(k+1)^r
\exp\left(-\frac{C}{2}(k+1)^\beta\right)\exp\left(-\frac{C}{2}(k+1)^\beta\right).
$$
An undergraduate analysis of the function
$u(x)=x^r\exp(-\frac{C}{2}x^\beta)$ gives
$$\max_{x\geq 0 }u(x)=u(x_0)
=\left(\frac{2r}{C\beta}\right)^{r/\beta}e^{-r/\beta}=
\left(\frac{2}{C\beta e}\right)^{r/\beta} r ^{r/\beta},\quad
x_0=\left(\frac{2r}{C\beta}\right)^{1/\beta},
$$
so that
$$ M_r\leq B\left(\frac{2}{C\beta e}\right)^{r/\beta}
r^{r/\beta},\qquad
B=C\sum_{k=0}^\infty\exp\left\{-\frac{C}{2}(k+1)^\beta\right\}.
$$
Hence,
$$ M_{k_1+1} \ldots M_{k_p+1} \leq  B^p \left(\frac{2}{C\beta
e}\right)^{\frac{n+p}{\beta}} (k_1+1)^{\frac{k_1+1}{\beta}} \ldots
(k_p+1)^{\frac{k_p+1}{\beta}}\,. $$ The inequalities
$$ \left(\frac{n+1}{n}\right)^{n/\beta}<e^{1/\beta},\qquad
(n+1)^{1/\beta}<e^{(n+1)/\beta},\qquad
(n+1)^{(n+1)/\beta}<e^{1/\beta}n^{n/\beta}e^{(n+1)/\beta}
$$
lead to the bound
$$ M_{k_1+1} \ldots M_{k_p+1} \leq B_2
\left(\frac{2}{C\beta}\right)^{n/\beta}k_1^{k_1/\beta} \ldots
k_p^{k_p/\beta}\,.
$$
Substituting the latter into (\ref{3.10}) gives
\begin{equation}
\max_{z\in\overline{\mathbb{D}}} |\gamma^{(n)}(z)| \leq C 4^n p^n
\left(\frac{2}{C\beta}\right)^{n/\beta} n^{k_1/\beta} \ldots
n^{k_p/\beta} \leq C \tilde C^n n^{n/\beta}, \quad n\ge 0.
\end{equation}
In other words, the function $\gamma$ for $D\in\mathcal{P}(\beta)$
belongs to the {\it Gevr\'e class} $\mathcal{G}_\beta$.

The rest of the proof goes along the same line of reasoning as in
\cite{EG2}. Suppose first that $0<\beta<1/2$. The celebrated Carleson's
theorem \cite{C} gives a complete description of zero sets for
$\mathcal{G}_\beta$. It claims that
$$\sum_{j=1}^\infty |l_j|^\frac{1-2\beta}{1-\beta} <\infty,$$
where $\{l_j\}$ are adjacent arcs of the zeros set of $\gamma$.
Hence the right inequality in (\ref{3.8z}) is obtained. The left
inequality is the general fact of fractal dimension theory.

When $\beta=1/2$ the Gevr\'e class $\mathcal{G}_{1/2}$ is known to
be quasi-analytic, i.e., it contains no nontrivial function $f$,
such that $f^{(n)}(\zeta_0)=0$ for all $n\geq0$ and some
$\zeta_0\in\mathbb{T}$. So each function $f\in\mathcal{G}_{1/2}$
may have only finite number of zeros inside the unit disk which
due to the relation between the discrete spectrum of $D$ and zeros
of $\gamma$ proves the second statement of the theorem. \hfill
$\square$

It turns out that the exponent $1/2$ in theorem 1 is sharp in the
following sense.
\begin{theorem}
\label{theor4} For arbitrary $\varepsilon > 0$ and different
points $-2<\nu_i<2, \ i=1,\ldots,p,$ there exists an operator
$D\in \mathcal{P}_p(\frac{1}{2} - \varepsilon)$ such that its
discrete spectrum $\sigma_{\rm d}(D)$ is infinite and, moreover,
$E_D=\{\nu_1,\ldots,\nu_p\}$.
\end{theorem}
{\it Proof.}  The proof of the theorem is based on the similar
result, proved in \cite{EG2} for the Jacobi matrices, that is, for $p=1$.
If $p>1$ than given $\nu_1,\ldots,\nu_p$ we construct the complex
Jacobi matrices
\begin{equation}
   J(i)=\left(\begin{array}{ccccc}
                b_0(i)&a_0(i)&&&\\
                a_0(i)&b_1(i)&a_1(i)&&\\
                &a_1(i)&b_2(i)&a_2(i)&\\
                &&\ddots&\ddots&\ddots\\
            \end{array}\right),
   \quad \begin{array}{c} a_j(i)>0,\, b_j(i) \in \mathbb{R},\\
                          j\ge 1,\ i=1,2,\ldots,p\\
   \end{array}
\end{equation}
which belong to the class $\mathcal{P}_1(\frac{1}{2} -
\varepsilon)$ and have the only point of accumulation of the
discrete spectrum $E_{J(i)}=\nu_i$.

Consider now the $p$-banded matrix $D=\|d_{ij}\|_{i,j=1}^\infty$
with
$$ d_{pn+i, \, pn+i} = b_n(i), \quad
   d_{pn+i, \, p(n+1)+i} = d_{p(n+1)+i, \, pn+i} = a_n(i),
$$
$$ d_{pn+i, \, pn+i+j} = d_{pn+i+j, \, pn+i} = 0, \qquad j\ne 0,p,
   \quad i=1,\ldots,p.
$$
Since $J(i)\in\mathcal{P}_1(\frac{1}{2} - \varepsilon)$, it is
easily seen that $D \in \mathcal{P}_p(\frac{1}{2} - \varepsilon)$.

The space $\ell^2$, where the operator $D$ acts, can be decomposed
as follows:
$$ \ell^2 = \bigoplus_{i=1}^p L_i
$$
with $ L_i={\rm Lin} \{\hat e_i, \hat e_{p+i}, \hat e_{2p+i},\ldots,
\hat e_{np+i},\ldots \}$ and $\{\hat e_k\}_{k\in\mathbb{N}}$ the
standard basis vectors in $\ell^2$. It is shown directly that the
subspaces $L_i$ are invariant for $D$.

The restriction $D_i=D\left|L_i\right.$ of $D$ on the subspace
$L_j$ has the matrix representation $J(i)$. Since
$D=\bigoplus_{i=1}^p D_i$ and every operator $D_i$ has discrete
spectrum with the only accumulation point $\nu_i$, the discrete
spectrum $D$ is the union of discrete spectrum of the operators
$J(i)$ and it has $p$ accumulation points $\nu_1, \nu_2, \ldots,
\nu_p$, as needed. \hfill $\square$

\medskip

The second issue we address in this section concerns domains which
contain the whole discrete spectrum, and conditions for the lack
of the discrete spectrum.

We start out with (\ref{2.24x}) and (\ref{2.26a}) for $n=0$
\begin{equation} \label{3.1a}
\|\Delta(z)-I\|\le \frac{2\sigma_0(0)}{|z-z^{-1}|}\exp
\left\{\frac{2\sigma_0(0)}{|z-z^{-1}|}\right\}, \ \ \|\Delta
(z)-I\|\le \sigma_1(0)\exp\{\sigma_1(0)\}\,,
\end{equation}
with $\Delta=\tilde V_0$, $\sigma_0(0)=\sum_{k=1}^\infty h_k$,
$\sigma_1(0)=\sum_{k=1}^\infty kh_k$, $h_k$ defined in
(\ref{2.16}), which hold under assumptions (\ref{2.17}) and
(\ref{2.19}), respectively. To work with (\ref{3.1a}) we have to
find the {\it efficient} constant $C$ which enters (\ref{2.16.1}).
To that end let us go back to (\ref{2.5x}) and assume that
\begin{equation} \label{effic}
q:=\sup_{n\ge 1} q_n=\sup_{n\ge1} \Bigl(\left|d_{n,n-p} -1\right|
+ \left|d_{n,n+p}- 1\right| + \sum_{r={-p+1}}^{p-1} \left|
d_{n,n+r} \right|\Bigr)<1,
\end{equation}
which now implies
\begin{equation} \label{2.26x}
\sup_{k\ge 1}\|A_k-I\|\le q, \qquad \sup_{k\ge 1}\|A_k\|\le 1+q<2.
\end{equation}
It is not hard to show now that $\sup_{k\ge 1}\|A_k^{-1}\|\le
(1-q)^{-1}$,
\[
\sup_{k\ge 1}\|L_k\|\le \exp\{\sum_{j=1}^\infty \|A_j-I\|\},
\qquad \sup_{k\ge 1}\|L_k^{-1}\|\le
\exp\Bigl\{\frac1{1-q}\,\sum_{j=1}^\infty \|A_j-I\|\Bigr\}
\]
and
\[
\sup_{k\ge 1}\|L_k\|\|L_k^{-1}\|\le
\exp\Bigl\{\frac{2-q}{1-q}\,\sum_{j=1}^\infty \|A_j-I\|\Bigr\}\le
\exp\Bigl\{\frac{2-q}{1-q}\,Q_0\Bigr\}, \ \ Q_0:=\sum_{k=1}^\infty
q_k\,.
\]
Hence, instead of (\ref{2.16.1}), we have by (\ref{2.16x})
\[
h_k\le 2\exp\Bigl\{\frac{2-q}{1-q}\,Q_0\Bigr\} (\hat q_{k-1}+\hat
q_k)
\]
and so
\begin{equation} \label{3.1b}
\sigma_0(0)\le 4\exp\Bigl\{\frac{2-q}{1-q}\,Q_0\Bigr\}\,Q_0, \quad
\sigma_1(0)\le 4\exp\Bigl\{\frac{2-q}{1-q}\,Q_0\Bigr\}\,Q_1, \ \
Q_1:= \sum_{k=1}^\infty kq_k.
\end{equation}

Let now $t$ be a unique root of the equation
\begin{equation}
te^t=1,\quad t\approx 0.567. \label{3.4}
\end{equation}
\begin{theorem}
Assume that $Q_0<\infty$ and $q<1$. Then the domain
$$ G(D)=\{z+z^{-1}:\ \ z\in\Omega\}
$$
with
$$\Omega:=\{z\in\mathbb D: |z-z^{-1}|>\frac{8Q_0}{t}\,
  \exp\Bigl\{\frac{2-q}{1-q}\,Q_0\Bigr\} $$
is free from the discrete spectrum $\sigma_{\rm d}(D)$. Moreover,
assume that $Q_1<\infty$. Then $D$ has no discrete spectrum as
long as
$$ \exp\Bigl\{\frac{2-q}{1-q}\,Q_0\Bigr\}\,Q_1<\frac4{t}. $$
\end{theorem}
{\it Proof}. It is clear from the first inequality in (\ref{3.1a})
that $\|\Delta(z)-I\|<1$ (and the more so, $\det\Delta(z)\ne0$)
whenever $|z-z^{-1}|>2\sigma_0(0)t^{-1}$. The rest follows
immediately from the first bound in (\ref{3.1b}) and the relation
between eigenvalues of $D$ and zeros of $\Delta$. The second
statement is proved in exactly the same way by using the second
inequality in (\ref{3.1a}) and the second bound in (\ref{3.1b}).
\hfill $\square$

\medskip

{\sl Remark}. Suppose that
$$ c=\frac{8Q_0}{t}\,\exp\Bigl\{\frac{2-q}{1-q}\,Q_0\Bigr\}<2. $$
 Then $\sigma_{\rm d}(D)$ is contained in the union of two symmetric
rectangles
$$\sigma_{\rm d}(D)\subset\left\{w:\,\sqrt{4-c^2}<|\mbox{Re}\,w|<\sqrt{4+c^2},\
\ |\mbox{Im}\,w|<\frac{c^2}{4}\right\}. $$

\setcounter{section}{4} \setcounter{equation}{0}
\section {The discrete spectrum of doubly-infinite banded matrices}

A doubly-infinite complex matrix
$D=\|d_{ij}\|_{i,j=-\infty}^\infty$ is called the {\it banded
matrix of order $p$} if
\begin{equation}
d_{ij}=0, \quad |i-j|>p, \qquad d_{ij}\ne 0, \quad |i-j|=p, \qquad
d_{ij} \in \mathbb{C}.
\end{equation}

There is a nice ``doubling method'' that relates a doubly-infinite
$p$-banded matrix to a certain semi-infinite banded matrix of
order $2p$. Indeed, let $\{e_k\}_{k\in\mathbb{Z}}$ be the standard
basis in $\ell^2(\mathbb{Z})$. Consider a transformation $U:
\ell^2(\mathbb{Z}) \longrightarrow \ell^2(\mathbb N)$, defined by
$$ U e_k = \hat e_{-2k}, \quad k<0,
   \qquad U e_k = \hat e_{2k+1}, \quad k\geq 0,
$$
where $\{\hat e_k\}_{k\in\mathbb N}$ is the standard basis in
$\ell^2(\mathbb N)$. The transformation $U$ is clearly isometric
and it is easy to check directly that the matrix
\begin{equation}\label{4.1x}
\hat D := U D U^{-1}
\end{equation}
is the semi-infinite band matrix of order $2p$ unitarily
equivalent to $D$. For instance, in the case of Jacobi matrices
($p=1$)
\[
   J=\left(\begin{array}{ccccccc}
            \ddots&\ddots&\ddots&&&\\
            &a_{-1}&b_{-1}&c_0&&\\
            &&a_0&b_0&c_1&\\
            &&&a_1&b_1&c_2&\\
            &&&&\ddots&\ddots&\ddots\\
            \end{array}\right),\ \
  \hat J=\left(\begin{array}{ccccc}
            \hat B_0& \hat C_1&&&\\
            \hat A_1& \hat B_1& \hat C_2&&\\
            & \hat A_2& \hat B_2& \hat C_3&\\
            &&\ddots&\ddots&\ddots\\
            \end{array}\right)
\]
with
\[
\begin{aligned}
\hat B_0 &=\left(\begin{array}{cc}
            b_0&a_0\\
            c_0&b_{-1}\\
            \end{array}\right), \quad
\hat B_k=\left(\begin{array}{cc}
            b_k&0\\
            0&b_{-k-1}\\
            \end{array}\right), \\
\hat A_k &=\left(\begin{array}{cc}
            a_k&0\\
            0&c_{-k}\\
            \end{array}\right), \quad
\hat C_k=\left(\begin{array}{cc}
            c_k&0\\
            0&a_{-k}\\
            \end{array}\right)  \quad k\in\mathbb N.
\end{aligned}
\]

Similarly to the semi-infinite case (\ref{0.3}), we put
\begin{equation}
q_n:= |d_{n, n-p}-1| + \sum_{j=-p+1}^{p-1} |d_{n, n+j}|+|d_{n,
n+p}-1|, \quad m\in \mathbb{Z}.
\end{equation}
We say that the matrix $D$ belongs to the class
$\mathcal{P}_p(\beta)$, $0<\beta<1$, if
\begin{equation}
q_n \leq C_1\,\exp(-C_2\, |n|^\beta), \qquad C_1, C_2>0, \qquad
n\in \mathbb{Z}.
\end{equation}
It is clear by the construction that the semi-infinite matrix
$\hat D$ (\ref{4.1x}) belongs to $\mathcal{P}_{2p}(\beta)$
whenever the doubly-infinite matrix $D$ belongs to
$\mathcal{P}_p(\beta)$. So the results of the previous section can
be easily derived for doubly-infinite banded matrices as well. For
instance, the following statement holds.

\begin{theorem} \label{theor5.1} Let
$D\in\mathcal{P}_p(\beta)$ where $0<\beta<\frac{1}{2}$. Then $E_D$
is a closed point set of the Lebesgue measure zero and its
convergence exponent satisfies
\begin{equation} \label{1.6} \dim E_D \leq
\tau(E_D)\leq\frac{1-2\beta}{1-\beta}\,,
\end{equation}
where $\dim E_D$ is the Hausdorff dimension of $E_D$. Moreover, if
$D\in \mathcal{P}_p(\frac{1}{2})$ then $E_D=\emptyset$ i.e. the
discrete spectrum is finite.
\end{theorem}

\setcounter{section}{5} \setcounter{equation}{0}
\section {The discrete spectrum of asymptotically periodic Jacobi matrices}

The goal of this section is to apply the results obtained above to
the study of the spectrum of doubly-infinite asymptotically
periodic complex Jacobi matrices.

\begin{definition}
A doubly-infinite complex Jacobi matrix
\begin{equation}\label{6.1}
   J^0=\left(\begin{array}{ccccccc}
            \ddots&\ddots&\ddots&&&\\
            &a_{-1}^0&b_{-1}^0&c_0^0&&\\
            &&a_0^0&b_0^0&c_1^0&\\
            &&&a_1^0&b_1^0&c_2^0&\\
            &&&&\ddots&\ddots&\ddots\\
            \end{array}\right), \quad b_0^0=(J^0e_0,e_0),
\end{equation}
$a_n^0c_n^0\not=0$, $a_n^0,b_n^0,c_n^0\in\mathbb C$ is called
$p$-periodic if
$$ a_{n+p}^0=a_n^0, \quad b_{n+p}^0=b_n^0, \quad
   c_{n+p}^0=c_n^0, \quad n\in\mathbb Z.
$$
\end{definition}

The following result due to P.B. Naiman \cite{N2,N} is the key
ingredient in our argument.

{\bf Theorem}. {\it Each $p$-periodic Jacobi matrix $J^0$
satisfies the algebraic equation of $p$-th degree
$$ P(J^0) = E_p,
$$
where $E_p = \{e_{jk}\}$ is a $p$-banded matrix with the entries
\begin{equation} \label{6.2}
 e_{ij} = \left\{ \begin{array}{ll}
     \alpha := {\displaystyle \prod_{k=1}^p a_k^0},& \quad i-j=p,\\
     \delta := {\displaystyle \prod_{k=1}^p c_k^0}, &\quad i-j=-p,\\
     0, &\quad |i-j| \ne p.\\
  \end{array}\right.
\end{equation}
Next, let
\[
\Pi:=\{z\in\mathbb C: z=\alpha e^{ipt}+\delta e^{-ipt}, \quad 0\le
t\le \frac{2\pi}p \}
\]
be an ellipse in the complex plane, and
\[
\Gamma:=P^{(-1)}(\Pi)=\{z\in\mathbb C: P(z)\in\Pi\}
\]
its preimage. Then the spectrum} $\sigma(J^0)=\Gamma$.

The polynomial $P$ is known as the {\it Burchnall--Chaundy
polynomial for $J^0$} and given explicitly in \cite[p.141] {N}:
\[
P(\lambda)=\det\bigl(\lambda-J(1,p)\bigr)-a_1c_1\det\bigl(\lambda-J(2,p-1)\bigr),
\ \
   J(m,n):=\left(\begin{array}{cccc}
            b_m^0&c_{m+1}^0&&\\
            a_{m+1}^0&b_{m+1}^0&\ddots&\\
            &\ddots&\ddots&c_n^0\\
            &&a_n^0&b_n^0\\
            \end{array}\right).
\]

{\sl Remark.} There is yet another way to exhibit the spectrum of
any $p$-periodic Jacobi matrix $J^0$. Consider the $p\times
p$-matrix
\[
   \Lambda_p(t):=\left(\begin{array}{cccc}
            b_1^0&c_2^0&&a_1^0e^{ipt}\\
            a_2^0&b_2^0&\ddots&\\
            &\ddots&\ddots&c_p^0\\
            c_1^0e^{-ipt}&&a_p^0&b_p^0\\
            \end{array}\right).
\]
with the eigenvalues $\{\lambda_j(t)\}_{j=1}^p$, and put
\[
{\rm range} \,\lambda_j:=\{w\in\mathbb C: w=\lambda_j(t), \ \ 0\le
t\le\frac{2\pi}p\}.
\]
Then
\[
\sigma(J^0)=\bigcup_{j=1}^p {\rm range}\,\lambda_j.
\]
Indeed, it is easy to compute
\[
\det\left(\Lambda(t)-\lambda\right)=(-1)^p\left(P(\lambda)-\alpha
e^{ipt}-\delta e^{-ipt}\right).
\]

For further information about the spectra of periodic Jacobi
matrices with algebro-geometric potential see \cite{BG}.

For the rest of the paper we restrict ourselves with the so called
{\it quasi-symmetric} matrices with
\begin{equation}\label{6.3}
 \prod_{k=1}^p a_k^0 = \prod_{k=1}^p c_k^0, \qquad \alpha=\delta.
\end{equation}
The latter certainly holds for the symmetric periodic matrices ($
a_n^0 =c_n^0$). In this case we have $Q(J^0)=D_0$ with $D_0$ defined
in (\ref{0.11}) for the polynomial $Q=\alpha^{-1}P$, and
\begin{equation}\label{6.3a}
\sigma(J^0)=Q^{(-1)}\bigl([-2,2]\bigr).
\end{equation}
It is not hard to make sure that $\sigma(J^0)$ is a collection of
finitely many algebraic arcs with no closed loops, that is, the
complement $\mathbb C\backslash\sigma(J^0)$ is a connected set.
Indeed, if some of these arcs formed a closed loop
$\Gamma_1\subset\Gamma$ with the interior domain $G_1$, then $\Im
Q$ would vanish on $\Gamma_1$ and so $\Im Q\equiv0$ in $G_1$,
which is impossible, unless $Q$ is a real constant.

Let now $J$ be a complex asymptotically $p$-periodic
doubly-infinite Jacobi matrix with the quasi-symmetric
$p$-periodic background $J^0$ (\ref{6.1})
\begin{equation}\label{6.4}
  J=\left(\begin{array}{ccccccc}
            \ddots&\ddots&\ddots&&&\\
            &a_{-1}&b_{-1}&c_0&&\\
            &&a_0&b_0&c_1&\\
            &&&a_1&b_1&c_2&\\
            &&&&\ddots&\ddots&\ddots\\
            \end{array}\right), \quad a_nc_n\not=0,
\end{equation}
and
$$ \omega_n:=|a_n-a_n^0| + |b_n-b_n^0| + |c_{n+1}-c_{n+1}^0| \to 0, \quad
k\to \pm\infty.
$$
In other words, $J-J^0$ is a compact operator in $\ell^2(\mathbb
Z)$. Since $\mathbb C\backslash\sigma(J^0)$ is connected, the
version of Weyl's theorem holds (cf. \cite{GK}, Lemma I.5.2)
\[
\sigma(J)=\sigma(J^0)\cup\sigma_{\rm d}(J),
\]
where the discrete spectrum $\sigma_{\rm d}(J)$ is an at most
denumerable set of eigenvalues of finite algebraic multiplicity
off $\sigma(J^0)$, which can accumulate only to $\sigma(J^0)$.

We proceed with two simple propositions.

\begin{proposition}\label{lem6.2}
Let $D_{m_1}$ and $D_{m_2}$ be banded matrices of orders $m_1$ and
$m_2$, respectively. Then
\begin{enumerate}
    \item  $D=D_{m_1}D_{m_2}$ is the banded matrix of order
    $m=m_1+m_2$.
    \item  If $m_1<m_2$ then $D_{m_1}+D_{m_2}$ is the banded
    matrix of order $m_2$.
\end{enumerate}
In particular, if $T$ is a polynomial of degree $p$ and $D$ is a
banded matrix of order $m$, then $T(D)$ is the banded matrix of
order $pm$.
\end{proposition}
{\it Proof.} (2) is obvious. To prove (1) let us show
that the elements of the extreme diagonals of $D$ do not vanish.
Indeed, let $D_{m_l}=\{d_{ij}^{(l)}\}$, $l=1,2$, and
$D=\{d_{ij}\}$. Then
$$ d_{k, \ k+m}= \sum_{j=1}^\infty d_{k\, j}^{(1)} d_{j,\,
k+m}^{(2)} = d_{k\, k+m_1}^{(1)} d_{k+m_1,\, k+m}^{(2)} \ne 0,
\quad k\in\mathbb N.
$$
Similarly, $d_{k,k-m}=d_{k,k-m_1}^{(1)}d_{k-m_1,k-m}^{(2)}\ne 0$.
The rest is plain. \hfill $\square$

\begin{proposition}\label{lem6.3}
Let $J_1$ and $J_2$ be bounded (with bounded entries) complex
Jacobi matrices. Put
\[
\omega_k=\sum_j \left|(J_1-J_2)_{kj}\right|=\sum_j\left| \langle
(J_1-J_2)e_j,e_k\rangle\right|.
\]
Then for any polynomial $T$ of degree $p$ we have
\[
\sum_j \left|(T(J_1)-T(J_2))_{kj}\right|=\sum_j\left| \langle
(T(J_1)-T(J_2))e_j,e_k\rangle\right|\leq
C\sum_{s=k-p}^{k+p}\omega_s,
\]
where a positive constant $C$ depends on $J_1$, $J_2$ and $T$.
\end{proposition}
{\it Proof.} For $T(z)=z^n$ the statement follows immediately from
the banded structure of the powers $J_l^m$, $(J_l^*)^s$, $l=1,2$
and the equality
\[
J_1^n-J_2^n=\sum_{j=0}^{n-1} J_1^{n-j-1}(J_1-J_2)J_2^j.
\]
The rest is straightforward. \hfill $\square$

\begin{definition}
We say that the matrix $J$ $(\ref{6.4})$ belongs to the class
$\mathcal P(\beta, J^0)$, if
\begin{equation}
\omega_n \leq C_1\,\exp(-C_2\, |n|^\beta), \qquad 0<\beta<1 ,
\quad C_1, C_2>0; \ n\in \mathbb{Z}.
\end{equation}
\end{definition}

Our main result claims that $\sigma_{\rm d}(J)$ is a finite set as
long as $J\in\mathcal P(1/2, J^0)$.

\begin{theorem}
Let $J$ be an asymptotically $p$-periodic doubly-infinite Jacobi
matrix $(\ref{6.4})$ with the quasi-symmetric background $J^0$
$(\ref{6.1})$, $(\ref{6.3})$. If $J\in\mathcal P(1/2, J^0)$, then
$\sigma(J)$ is the union of $p$ algebraic arcs and a finite number
of eigenvalues of finite algebraic multiplicities off these arcs.
\end{theorem}
{\it Proof.}  Let $T=Q=\alpha^{-1}P$ be the Burchnall--Chaundy
polynomial for $J^0$ taken from (\ref{6.3a}). By Proposition
\ref{lem6.2} $Q(J)$ is the p-banded matrix. Since $J\in\mathcal
P(1/2, J^0)$ the matrix $Q(J)$ is close to $Q(J^0)=D_0$
with $D_0$ defined in (\ref{0.11}) in the sense of Proposition \ref{lem6.3}
\[
\sum_j \left|(Q(J)-D_0)_{kj}\right|\le C_1 e^{-C_2 |k|^{1/2}}
\]
and so $Q(J)\in\mathcal{P}_p(1/2)$. According to Theorem
\ref{theor5.1} $\sigma(Q(J))=[-2,2]\cup E$, where the set $E$ of
eigenvalues off $[-2,2]$ is now finite.

On the other hand, as we know, $\sigma(J)=\Gamma\cup F$, where the
set $F$ of eigenvalues off $\Gamma=\sigma(J^0)$ is at most
denumerable, and $\Gamma=Q^{(-1)}([-2,2])$. Therefore, by the
Spectral Mapping Theorem $Q(F)=E$ and so $F$ is a finite set, as
claimed. \hfill $\square$

\medskip

{\bf Acknowledgements.}  The authors thank Prof. M.Pituk for his
helpful comments on the asymptotic behavior of the solutions of
linear difference equations and for the reference to \cite{Perron}.


\begin{thebibliography}{}



\bibitem{EG1} {\it I.Egorova and L.Golinskii,} On location of
discrete spectrum for complex Jacobi matrices. ---
 {\it Proceedings of AMS}, {\bf 133} (2005), ¹ 12, 3635-3641.

\bibitem{EG2} {\it I.Egorova and L.Golinskii,} On limit sets for the discrete spectrum of
complex Jacobi matrices --- {\it Matematicheskiy sbornik}. {\bf
196} (2005), ¹ 6, 43--70.

\bibitem{P1} {\it B.S.Pavlov,} On nonselfadjoint Schr\"odinger operator I. ---
{\it Problems of mathematical physics}, LGU, {\bf 1} (1966),
102--132 .

\bibitem{P2} {\it B.S.Pavlov,} On nonselfadjoint Schr\"odinger
operator II, --- {\it Problems of mathematical physics}, LGU, {\bf
2} (1967), 135--157 .

\bibitem{GK} {\it I.Gohberg and M.Krein,} Introduction to the theory
of linear nonselfadjoint operators in Hilbert space. ``Nauka'',
Moscow, 1965 (Russian).

\bibitem{Perron} {\it O.Perron,} \"{U}ber Summengleichungen und
Poincar\'{e}sche Differenzengleichingen. --- {\it Math. Annalen},
{\bf 84} (1921), 1--15 (German).

\bibitem{C} {\it L.Carleson,} Sets of uniqueness for functions
analytic in the unit disc. --- {\it Acta Math}, {\bf 87} (1952),
325--345.

\bibitem{N2} {\it P.B.Naiman,} On the theory of periodic and limit-periodic Jacobi operators.
--- {\it Sov. Math. Dokl}, {\bf 3} (1962), 4, 383--385 (Russian).

\bibitem{N} {\it P.B.Naiman,} To the spectral theory of the
non-symmetric periodic Jacobi matrices. --- {\it Notes of the
Department for Mathematics and Mechanics of the Kharkov's state
University and of the Kharkov's mathematical Society}, {\bf XXX}
(1964), 4, 138--151 (Russian).

\bibitem{BG} {\it V. Batchenko and F. Gesztesy.} On the spectrum
of Jacobi operators with quasi-periodic algebro-geometric
coefficients, to appear in IMRN.

\end{thebibliography}
\end{document}